\documentclass[a4paper]{article}
\usepackage{amsmath,amssymb,amsfonts}

\newtheorem{theorem}{Theorem}[section]

\title{A proof of the non existence of Frey curves without using TSW theorem}
\author{Jailton C. Ferreira}

\date{ }
\begin{document}
\maketitle
\pagenumbering{arabic}

\begin{abstract}
Fermat's Last Theorem (FLT) implies that the Frey curves do not
exist. A proof of FLT independent of proved Taniyama-Shimura-Weil
conjecture is presented.
\end{abstract}

\section{Introduction} \label{sec-1}

\begin{quote}
\textit{Cubum autem in duos cubos, aut quadratoquadratum in duos
quadratoquadratos, et generaliter nullam in infinitum ultra
quadratum potestatem in duos ejusdem nominis fas est dividere:
cujus rei demonstrationem mirabilem sane detexi. Hanc marginis
exiguitas non caperet.}
\end{quote}

This quote was written by Pierre de Fermat (1601 - 1665) in the
margin of his copy of Diophantus' Arithmetikae, near to a problem
on Pythagorean triples.

\hspace{22pt} The Latin text states that it is impossible to
separate a cube into two cubes, or a biquadrate into two
biquadrates, or in general any power higher than the second into
two powers of the like degree, in other words, for any integer $n
> 2$, the equation
\begin{equation}\label{um-1}
a^n + b^n = c^n
\end{equation}

does not have any integer solutions with $abc \ne 0$. The quote
ends with
\begin{quote}
I have discovered a truly marvelous proof of this, which this
margin is too narrow to contain.
\end{quote}

Fermat's personal comment was discovered after his death. Now it
is generally accepted that he likely found a mistake in the proof.
However, he enjoyed issue theorems without the proof as a
challenge problem to other mathematicians and it is possible that
Fermat had outlined a proof. The personal note became the
conjecture known as Fermat's Last Theorem (FLT). ``Last Theorem"
because it was for many years the last remaining statement in the
list of Fermat's works to be independently proven or disproven.

\hspace{22pt} During 324 years after the publication of the
Fermat's note in 1670, many mathematicians failed in obtain a
proof of FLT. However, part of modern number theory was developed
in these attempts and the conjecture was proved for some values of
$n$. Euler proved \eqref{um-1} for $n = 3$. Sophie Germain proved
that if $n$ and $p=2n+1$ are primes, the equation \eqref{um-1} has
no solution with $abc \ne 0$. Dirichlet and Legendre proved for $n
= 5$. Lam\'e proved for $n = 7$. Kummer found that FLT is true for
all regular prime exponents. Wieferich proved case I of FLT for
$n$, where $n$ does not divide the Fermat quotient $(2^{n-1} -
1)/n$. Vandiver proved FLT for certain irregular prime exponents.
Gerd Faltings ~\cite{Faltings} proved that \eqref{um-1} has at
most finitely many integers solutions. The conjecture was
validated with the help of computers up to 125000
~\cite{Ribenboim} and is true for $n < 4 \times 10^6$
~\cite{Buhler}.

\hspace{22pt} In 1984, Gerhard Frey proved  that for a prime $n
\ge 5$ the following three conditions are equivalent
~\cite{Frey1}:
\begin{enumerate}
\item There exists a solution of Fermat's equation \eqref{um-1}.
\item There exists a stable elliptic curve $E$ over $\mathbb{Q}$ such that
\begin{itemize}
\item[-] the points of order two of $E$ are $\mathbb{Q}$-rational,
\item[-] the field $K_n$ obtained by adjoining the coordinates
of the points of order $n$ to $\mathbb{Q}$ is unramified outside
$2 \cdot n$,
\item[-] Min$(0,v_{n}(j_{E})) \equiv 0 \equiv v_{2}(j_{E}) - 8$
mod $n$.
\end{itemize}
\item There exists a stable elliptic curve $E$ over $\mathbb{Q}$ with a
minimal equation
\begin{equation}\label{um-2}
y^2 + xy = x^3 \alpha x^2 + \beta x \hspace{8pt} \textrm{with}
\hspace{8pt} \alpha , \beta \in \mathbb{Z}
\end{equation}
and
\begin{equation}\label{um-3}
2^{8} \cdot \Delta_{E} \in \mathbb{Z}^{2n}
\end{equation}
\end{enumerate}

Frey presented reasons for believing in the truth of the
conjecture ``\textit{E is not modular}". He concluded that this
conjecture and the conjecture of Taniyama-Shimura-Weil (conjecture
TSW) would imply Fermat's Last Theorem. The conjecture TSW states
that ``\textit{every elliptic curve over $\mathbb{Q}$ is a modular
elliptic curve}".

\hspace{22pt} The conjecture that $E$ is not modular was nearly
proved by Serre who left one part unproved, (Serre's
$\epsilon$-conjecture). In 1986, Ken Ribet proved that the
elliptic curve $E$ is not modular ~\cite{Ribet}. In 1994, Andrew
Wiles, assisted by Richard Taylor, proved that all semistable
curves over $\mathbb{Q}$ are modular ~\cite{Wiles} and
~\cite{Wiles&Taylor}. The proof of Fermat's Last Theorem was
completed.

\section{The non existence of $E$}
\label{sec-2}

\hspace{22pt} The way indicated by Frey to proof of Fermat's Last
Theorem requires the proof that the elliptic curves $E$ do not
exist. An independent way to prove FLT implies, for example, in
another proof that Frey curves do not exist.

\begin{theorem} \label{teorema-1}
For any integer $n
> 2$, the equation
\begin{equation}\label{dois-1}
x^n + y^n = z^n
\end{equation}

does not have any positive integer solutions with $xyz \ne 0$.
\end{theorem}
\textit{Proof:}

\hspace{22pt} Let be $x$, $y$, $z$, $n \in \mathbb{Z}$ and $x >
0$, $y > 0$, $z > 0$ and $n > 2$.

\hspace{22pt} Let $n$ be equal to $pq$, where $p$ is prime.
Substituting $n$ by $pq$ into \eqref{dois-1} we have
\begin{equation}\label{dois-2}
{(x^q)}^p + {(y^q)}^p = {(z^q)}^p
\end{equation}

From \eqref{dois-2} we verified that if FLT is true for $n = p$,
then it is true for $n = pq$. The proof for $n=4$ includes the
cases with $n = 8, 12, 16, 20, \ldots $ We need to prove FLT for
odd primes $n$ only.

\hspace{22pt} For $x = y$ we obtain $z = 2^{\frac{1}{n}}x$, that
is, there is no integer solution.

\hspace{22pt} Let be
\begin{equation}\label{dois-3}
u^p + v^p = 1
\end{equation}

where $p$ is prime greater than 2 and $u$, $v \in \mathbb{R}$ with $0<u<1$ and $0<v<1$. Let
us consider the line
\begin{equation}\label{dois-4}
v = -du + 1
\end{equation}
where
\begin{equation}\label{dois-5}
0 < d < 1
\end{equation}

For given values of $p$ and $u$, only one
line satisfying \eqref{dois-4} intercept \eqref{dois-3}. There exists only one slope of \eqref{dois-4} satisfying this
condition.

\hspace{22pt} Substituting $v$ from \eqref{dois-4} into
\eqref{dois-3} we have
\begin{equation}\label{dois-6}
u^p + (-du + 1)^p = 1
\end{equation}

Using the binomial theorem we obtain
\begin{equation}\label{dois-7}
{\Big( -du + 1 \Big)}^p = \sum_{k=0}^p  {p \choose k} (-du)^k
\end{equation}

From \eqref{dois-6} and \eqref{dois-7} we have
\begin{equation}\label{dois-8}
u^p + 1 + {p \choose 1} (-du)^1 + {p \choose 2} (-du)^2 + \ldots +
{p \choose {p-1}} (-du)^{p-1} + (-du)^p = 1
\end{equation}
or
\begin{equation}\label{dois-9}
d^p - {p \choose {p-1}}u^{-1}d^{p-1} + {p \choose
{p-2}}u^{-2}d^{p-2} - \ldots - {p \choose 2}u^{2-p}d^2 + {p
\choose 1}u^{1-p}d - 1 = 0
\end{equation}

\hspace{22pt} Applying Descarte's Rule of Signs to \eqref{dois-9}
we verify that the number of variations in sign of the polynomial
equation is $p$, hence it has $p$ positive roots or less than $p$
by an even number. The number 0 is not a root of \eqref{dois-9}.
The negatives roots of \eqref{dois-9} are the positive roots of
the equation obtained substituting $d$ by $-d$ into
\eqref{dois-9}. Doing this and applying Descarte's Rule of Signs,
we verify that there are no negative roots. By Du Gua-Huat-Euler
theorem ~\cite{Mignotte} a polynomial equation
\begin{equation}\label{dois-10}
a_{p}d^p + a_{p-1}d^{p-1} + a_{p-2}d^{p-2} + \ldots + a_{2}d^2 +
a_{1}d + a_0 = 0
\end{equation}
has no complex roots if
\begin{equation}\label{dois-11}
(a_{k})^2 > a_{k-1}a_{k+1}
\end{equation}

for all $0 < k < p$. Substituting \eqref{dois-10} by
\eqref{dois-9} we obtain
\begin{equation}\label{dois-12}
{p \choose {k}}^2 > {p \choose
{k-1}}\times {p \choose {k+1}}
\end{equation}
 or
\begin{equation}\label{dois-13}
\Bigl[(p-k)! \times k!\Bigr]^2 < (p-k-1)! \times (k+1)! \times
(p-k+1)! \times (k-1)!
\end{equation}
or
\begin{equation}\label{dois-14}
(p-k)k < (p-k+1)(k+1)
\end{equation}

We can conclude from \eqref{dois-14} that \eqref{dois-9} has not complex roots.

\hspace{22pt} Let us consider the case where $u$ is $\frac
{x}{z}$. Let us notice that
\begin{equation}\label{dois-17}
0 < \frac {x}{z} < 1
\end{equation}

Substituting $u$ by $\frac {x}{z}$ into \eqref{dois-9} we have
\begin{equation}\label{dois-18}
d^p - {p \choose {p-1}}{\Big(\frac {x}{z}\Big)}^{-1}d^{p-1} +
\ldots - {p \choose 2}{\Big(\frac {x}{z}\Big)}^{2-p}d^2 + {p
\choose 1}{\Big(\frac {x}{z}\Big)}^{1-p}d - 1 = 0
\end{equation}

There exists only one slope $-d$ for this case, which is the slope
of the line that passes through the points (0, 1) and ($\frac
{x}{z}$, $\alpha$), where
\begin{equation}\label{dois-19}
\alpha = {\Big( 1 - {\Big(\frac {x}{z}\Big)}^p
\Big)}^{\frac{1}{p}}
\end{equation}

The slope is
\begin{equation}\label{dois-20}
\beta = \frac{(\alpha -1)z}{x}
\end{equation}

\hspace{22pt} For the polinomial equation \eqref{dois-10} part of the Vieta's formula is

\begin{equation}\label{dois-21}
a_0 = a_p (-1)^p \beta_1 \beta_2 ... \beta_p
\end{equation}
where the $\beta_{i}'s$ are the roots. To be consistent with only one value $\beta$, all
roots of the equation \eqref{dois-18} must be equal. Considering this and the values of $a_0$ and $a_p$ we obtain
\begin{equation}\label{dois-22}
(-1)^p (\beta)^p = -1
\end{equation}
since that $p$ is an odd number we have
\begin{equation}\label{dois-23}
\beta = 1
\end{equation}
However, \eqref{dois-23} does not satisfy the
condition \eqref{dois-5}. Therefore, does not exist solution for
\eqref{dois-18}.

\hspace{22pt} If FLT is not true there exists some rational number
$\beta$ equal to
\begin{equation}\label{dois-24}
\frac{\frac{y}{z} - 1}{\frac{x}{z}}
\end{equation}
where
\begin{equation}\label{dois-25}
\Big(\frac{y}{z}\Big)^p + \Big(\frac{x}{z}\Big)^p = 1
\end{equation}

that is solution for \eqref{dois-18}. Since does not exist
solution for \eqref{dois-18}, we conclude that for any integer $n
> 2$ the equation \eqref{dois-1} does not have any positive
integer solutions with $xyz \ne 0$.


\begin{thebibliography}{99}
\addcontentsline{toc}{chapter}{Bibliography}
\bibitem{Faltings}  G. Faltings, \textit{Endlichkeitss\"atze f\"ur abelsche Variet\"aten \"uber Zahlk\"orpern}, Acta
Inv. Math. \textbf{73} (1983), 349-366 (English translation in
~\cite{Cornell&Silverman}).
\bibitem{Cornell&Silverman}  G. Cornell and J. H. Silverman, eds.,
\textit{Arithmetic Geometry}, Springer-Verlag, New York, Berlin,
Heidelderg (1986).
\bibitem{Ribenboim}  P. Ribenboim, \textit{13 Lectures on Fermat's Last Theorem}, Springer-Verlag, New York
(1979).
\bibitem{Buhler}  J. Buhler, R. Crandall, R. Ernvall, T. Mets\"ankyl\"a, \textit{Irregular primes and ciclotomic invariants to four million}, Math. Comp. \textbf{61} (1993),
151-153.
\bibitem{Frey1}  G. Frey,
\textit{Links Between Stable Elliptic Curves And Certain
Diophantine Equations}, Ann. Univ. Sarav. Math. Ser. \textbf{1}
(1986), 1-40.
\bibitem{Ribet}  K. Ribet,
\textit{On modular representations of Gal($\overline{Q}$/Q)
arising from modular forms}, Inventiones Mathematikae \textbf{100}
(1990), 431-476.
\bibitem{Wiles}  A. Wiles,
\textit{Modular elliptic curves and Fermat's Last Theorem}, Annals
of Mathematics \textbf{142} (1995), 443-551.
\bibitem{Wiles&Taylor}  R. Taylor, A. Wiles,
\textit{Ring-theoretics properties of certain Hecke algebras},
Annals of Mathematics \textbf{142} (1995), 553-572.
\bibitem{Mignotte}  M. Mignotte,
\textit{Math\'ematiques pour le calcul formel}, Presses
Universitaires of France, Paris (1989).


\end{thebibliography}
\end{document}